\documentclass[11pt,a4paper]{amsart}
\usepackage{graphicx}
\usepackage{amssymb}
\usepackage{amsmath}
\usepackage{amsthm}
\usepackage{amscd}
\usepackage[all,2cell,dvips]{xy}
\UseAllTwocells \SilentMatrices
\newtheorem{thm}{Theorem}[section]

\newtheorem{lem}[thm]{Lemma}

\theoremstyle{definition}

\theoremstyle{remark}

\numberwithin{equation}{section}

\begin{document}
\title[An Auslander-type result for Gorenstein-projective modules]
{An Auslander-type result for Gorenstein-projective modules}
\author[  Xiao-Wu Chen
] {Xiao-Wu Chen}
\thanks{This project was supported by China Postdoctoral Science Foundation No. 20070420125, and
was also partially supported by the National Natural Science
Foundation of China (Grant No.s 10725104, 10501041 and 10601052).}
\thanks{E-mail:
xwchen$\symbol{64}$mail.ustc.edu.cn}
\maketitle
\date{}%
\dedicatory{}%
\commby{}%
\begin{center}
Department of Mathematics\\
 University of Science and
Technology of China \\Hefei 230026, P. R. China
\end{center}

\begin{abstract}
An artin algebra $A$ is said to be CM-finite if there are only
finitely many, up to isomorphisms, indecomposable finitely generated
Gorenstein-projective $A$-modules. We prove that for a Gorenstein
artin algebra, it is CM-finite if and only if every its
Gorenstein-projective module is a direct sum of finitely generated
Gorenstein-projective modules. This is an analogue of Auslander's
theorem on algebras of finite representation type (\cite{A,A1}).
\end{abstract}

\section{Introduction}

Let $A$ be an artin $R$-algebra, where $R$ is a commutative artinian
ring. Denote by $A\mbox{-Mod}$ (resp. $A\mbox{-mod}$) the category
of (resp. finitely generated) left $A$-modules. Denote by
$A\mbox{-Proj}$ (resp. $A\mbox{-proj}$) the category of (resp.
finitely generated) projective $A$-modules. Following \cite{IK}, a
chain complex $P^\bullet$ of projective $A$-modules  is defined to
be \emph{totally-acyclic}, if for every projective module $Q\in
A\mbox{-Proj}$ the Hom-complexes ${\rm Hom}_A(Q, P^\bullet)$ and
${\rm Hom}_A(P^\bullet, Q)$ are exact. A module $M$ is said to be
\emph{Gorenstein-projective} if there exists a totally-acyclic
complex $P^\bullet$ such that the $0$-th cocycle $Z^0(P^\bullet)=M$.
Denote by $A\mbox{-GProj}$ the full subcategory of
Gorenstein-projective modules. Similarly, we define finitely
generated Gorenstein-projective modules by replacing all modules
above by finitely generated ones, and we also get the category
$A\mbox{-Gproj}$ of finitely generated Gorenstein-projective modules
\cite{EJ}. It is known that $A\mbox{-Gproj}=A\mbox{-GProj}\cap
A\mbox{-mod}$ (\cite{C1}, Lemma 3.4). Finitely generated
Gorenstein-projective modules are also referred as maximal
Cohen-Macaulay modules. These modules play a central role in the
theory of singularity \cite{Buc, BEH, Bel3, C1} and of relative
homological algebra \cite{Bel1, EJ}.

\vskip 5pt

An artin algebra $A$ is said to be \emph{CM-finite} if there are
only finitely many, up to isomorphisms, indecomposable finitely
generated Gorenstein-projective modules. Recall that an artin
algebra $A$ is said to be of \emph{finite representation type} if
there are only finitely many isomorphism classes of indecomposable
finitely generated modules. Clearly, finite representation type
implies CM-finite. The converse is not true, in general. \vskip 5pt

 Let us recall the following famous result of Auslander \cite{A, A1}
 (see also Ringel-Tachikawa \cite{RT}, Corollary 4.4) : \vskip 10pt

 \noindent {\bf Auslander's
Theorem}\; An artin algebra $A$ is of finite representation type if
and only if every $A$-module is a direct sum of finitely generated
modules, that is, $A$ is left pure semisimple, see \cite{S1}.

\vskip 10pt

Inspired by the theorem above, one may conjecture the following
Auslander-type result for Gorenstein-projective modules: an artin
algebra $A$ is CM-finite if and only if every Gorenstein-projective
$A$-module is a direct sum of finitely generated ones. However we
can only prove this conjecture in a nice case.

\vskip 5pt

 Recall that an artin algebra $A$ is said to be Gorenstein
 \cite{Ha2}
if the regular module $A$ has finite injective dimension both at the
left and right sides. Our main result is

\vskip 10pt

\noindent {\bf Main Theorem }\; Let $A$ be a Gorenstein artin
algebra. Then $A$ is CM-finite if and only if every
Gorenstien-projective $A$-module is a direct sum of finitely
generated Gorenstein-projective modules.

\vskip 10pt

Note that our main result has a similar character to a result by
Beligiannis (\cite{Bel1}, Proposition 11.23), and also note that
similar concepts were introduced and then  similar results and ideas
were developed by Rump in a series of papers \cite{Ru1, Ru2, Ru3}.

\section{Proof of Main Theorem}

Before giving the proof, we  recall some notions and known results.

\subsection{} Let $A$ be an artin $R$-algebra. By a subcategory
$\mathcal{X}$ of $A\mbox{-mod}$, we mean a full additive subcategory
which is closed under taking direct summands. Let $M \in
A\mbox{-mod}$. We recall from \cite{AS1,AR1} that a \emph{right
$\mathcal{X}$-approximation} of $M$ is a morphism $f:
X\longrightarrow M$ such that $X\in \mathcal{X}$ and every morphism
from an object in $\mathcal{X}$ to $M$ factors through $f$. The
subcategory $\mathcal{X}$ is said to be \emph{contravariantly-finite
in $A\mbox{\rm -mod}$} if each finitely generated modules has a
right $\mathcal{X}$-approximation. Dually, one defines the notions
of \emph{left $\mathcal{X}$-approximations} and
\emph{covariantly-finite} subcategories. The subcategory
$\mathcal{X}$ is said to be \emph{functorially-finite in $A\mbox{\rm
-mod}$} if it is contravariantly-finite and covariantly-finite.
Recall that a morphism $f: X\longrightarrow M$ is said to be
\emph{right minimal}, if for each endomorphism $h: X \longrightarrow
X$ such that $f=f\circ h$, then $h$ is an isomorphism. A right
$\mathcal{X}$-approximation $f: X \longrightarrow M$ is said to be a
\emph{right minimal $\mathcal{X}$-approximation} if it is right
minimal. Note that if a right approximation exists, so does right
minimal ones; a right minimal approximation, if in existence, is
unique up to isomorphisms. For details, see \cite{AS1, AR1, ARS}.

\vskip 5pt

 The following fact is known.

\begin{lem}
Let $A$ be an artin algebra. Then \\
(1).\; The subcategory $A\mbox{\rm -Gproj}$ of $A\mbox{\rm -mod}$ is
closed under taking direct summands, kernels of
epimorphisms and extensions, and contains $A\mbox{\rm -proj}$.\\
(2). \;  The category $A\mbox{\rm -Gproj}$ is a Frobenius exact
category \cite{Ke3}, whose relative projective-injective objects are
precisely contained in $A\mbox{\rm -proj}$. Thus the stable category
$A\underline{\mbox{\rm -Gproj}}$ modulo projectives is a triangulated category.\\
(3). \; Let $A$ be Gorenstein. Then the subcategory $A\mbox{\rm
-Gproj}$ of $A\mbox{\rm -mod}$ is functorially-finite.\\
(4). \; Let $A$ be Gorenstein. Denote by $\{S_i\}_{i=1}^n$ a
complete list of pairwise nonisomorphic simple $A$-modules. Denote
by $f_i: X_i \longrightarrow S_i$ the right minimal $A\mbox{\rm
-Gproj}$-approximations. Then every finitely generated
Gorenstein-projective module $M$ is a direct summand of some module
$M'$, such that there exists a finite chain of submodules $0=M_0
\subseteq M_1 \subseteq \cdots \subseteq M_{m-1}\subseteq M_{m}=M'$
with each subquotient $M_j/M_{j-1}$ lying in $\{X_i\}_{i=1}^n$.
\end{lem}

\noindent{\bf Proof.}\quad Note that $A\mbox{-Gproj}$ is nothing but
$\mathcal{X}_\omega$ with $\omega=A\mbox{-proj}$ in \cite{AR1},
section 5. Thus (1) follows from \cite{AR1}, Proposition 5.1, and
(3) follows from \cite{AR1}, Corollary 5.10(1) (just note that in
this case, $_AA$ is a cotilting module).

  Since $A\mbox{-Gproj}$ is closed under extensions, thus it becomes an
exact category in the sense of \cite{Ke3}. The property of being
Frobenius and the characterization of projective-injective objects
follow directly from the definition, also see \cite{C1}, Proposition
3.1(1). Thus by \cite{Ha1}, chapter 1, section 2, the stable
category $A\underline{\mbox{\rm -Gproj}}$ is triangulated.

 By (1) and (3), we see that (4) is a special case of \cite{AR1}, Proposition 3.8. \hfill
$\blacksquare$

\vskip 10pt

 Let $R$ be a commutative artinian ring as above. An
additive category $\mathcal{C}$ is said be to \emph{$R$-linear} if
all its Hom-spaces are $R$-modules, and the composition maps are
$R$-bilinear. An $R$-linear category is said to be {\em hom-finite},
if all its Hom-spaces are finitely generated $R$-modules. Recall
that an \emph{$R$-variety} $\mathcal{C}$ means a hom-finite
$R$-linear category which is skeletally-small and idempotent-split
(that is, for each idempotent morphism $e: X \longrightarrow X$ in
$\mathcal{C}$, there exists $u: X \longrightarrow Y$ and $v: Y
\longrightarrow X$ such that $e=v\circ u$ and ${\rm Id}_Y=u\circ
v$). It is well-known that a skeletally-small $R$-linear category is
an $R$-variety if and only if it is hom-finite and Krull-Schmidt
(i.e., every object is a finite sum of indecomposable objects with
local endomorphism rings). See \cite{Rin}, p.52 or \cite{CYZ},
Appendix A. Then it follows that any factor category (\cite{ARS},
p.101) of an $R$-variety is still an $R$-variety.

\vskip 5pt

Let $\mathcal{C}$ be an $R$-variety. We will abbreviate the
Hom-space ${\rm Hom}_{\mathcal{C}}(X, Y)$ as $(X, Y)$. Denote by
$(\mathcal{C}^{\rm op}, R\mbox{-Mod})$ (resp. $(\mathcal{C}^{\rm
op}, R\mbox{-mod})$) the category of contravariant $R$-linear
functors from $\mathcal{C}$ to $R\mbox{-Mod}$ (resp.
$R\mbox{-mod}$). Then $(\mathcal{C}^{\rm op}, R\mbox{-Mod})$ is an
abelian category and $(\mathcal{C}^{\rm op}, R\mbox{-mod})$ is its
abelian subcategory. Denote by $(-, X)$ the representable functor
for each $X\in \mathcal{C}$. A functor $F$ is said to be
\emph{finitely generated} if there exists an epimorphism $(-, C)
\longrightarrow F$ for some object $C\in \mathcal{C}$; $F$ is said
to be \emph{finitely presented (= coherent)} \cite{Au2,A}, if there
exists an exact sequence of functors $(-, C_1) \longrightarrow (-,
C_0) \longrightarrow F \longrightarrow 0.$ Denote by ${\bf
fp}(\mathcal{C})$ the subcategory of $(\mathcal{C}^{\rm op},
R\mbox{-Mod})$ consisting of finitely presented functors. Clearly,
${\bf fp}(\mathcal{C})\subseteq (\mathcal{C}^{\rm op},
R\mbox{-mod})$. Recall the duality
\begin{align*}D={\rm Hom}_R(-,
E): R\mbox{-mod}\longrightarrow R\mbox{-mod},
\end{align*}
where $E$ is injective hull of $R/{\rm rad}(R)$ as an $R$-module.
Therefore, it induces duality $D: (\mathcal{C}^{\rm op},
R\mbox{-mod}) \longrightarrow (\mathcal{C}, R\mbox{-mod})$ and $D:
(\mathcal{C}, R\mbox{-mod}) \longrightarrow (\mathcal{C}^{\rm op},
R\mbox{-mod})$. The $R$-variety $\mathcal{C}$ is called a
\emph{dualizing $R$-variety} \cite{AR2}, if this duality preserves
finitely presented functors.

\vskip 10pt

 The following observation is important.

\begin{lem}
Let $A$ be a Gorenstein artin $R$-algebra. Then the stable category
$A\underline{\mbox{\rm -Gproj}}$ is a dualizing $R$-variety.
\end{lem}

\noindent {\bf Proof.} \quad Since $A\mbox{-Gproj}\subseteq
A\mbox{-mod}$ is closed under taking direct summands, thus
idempotents-split. Therefore, we infer that $A\mbox{\rm -Gproj}$ is
an $R$-variety, and its stable category $A\underline{\mbox{\rm
-Gproj}}$ is also an $R$-variety. By Lemma 2.1(3), the subcategory
$A\mbox{-Gproj}$ is functorially-finite in $A\mbox{-mod}$, then by a
result of Auslander-Smal{\o} (\cite{AS1}, Theorem 2.4(b))
$A\mbox{-Gproj}$ has almost-split sequences, and thus theses
sequences induce Auslander-Reiten triangles in
$A\underline{\mbox{-Gproj}}$ (Let us remark that it is Happel
(\cite{Ha2}, 4.7) who realized this fact for the first time). Hence
the triangulated category $A\underline{\mbox{-Gproj}}$ has
Auslander-Reiten triangles, and by a theorem of Reiten-Van den Bergh
(\cite{RV}, Theorem I.2.4) we infer that
$A\underline{\mbox{-Gproj}}$ has Serre duality. Now by \cite{IY},
Proposition 2.11 (or \cite{C2}, Corollary 2.6), we deduce that
$A\underline{\mbox{-Gproj}}$ is a dualizing $R$-variety. Let us
remark that the last two cited results are given in the case where
$R$ is a field, however one just notes that the results can be
extended to the case where $R$ is a commutative artinian ring
without any difficulty. \hfill $\blacksquare$

\vskip 5pt

For the next result, we recall more notions on functors over
varieties. Let $\mathcal{C}$ be an $R$-variety and let $F\in
(\mathcal{C}^{\rm op}, R\mbox{-Mod})$ be a functor. Denote by ${\rm
ind}(\mathcal{C})$ the complete set of pairwise nonisomorphic
indecomposable objects in $\mathcal{C}$. The \emph{support} of $F$
is defined to ${\rm supp}(F)=\{C\in {\rm ind}(\mathcal{C})\; |\;
F(C)\neq 0\}$. The functor $F$ is \emph{simple} if it has no nonzero
proper subfunctors, and $F$ \emph{has finite length} if it is a
finite iterated extension of simple functors. Observe that $F$ has
finite length if and only if $F$ lies in $(\mathcal{C}^{\rm op},
R\mbox{-mod})$ and ${\rm supp}(F)$ is a finite set. The functor $F$
is said to be \emph{noetherian}, if its every subfunctor is finitely
generated. It is a good exercise to show that a functor is
noetherian if and only if every ascending chain of subfunctors in
$F$ becomes stable after finite steps (one may use the fact: for a
finitely generated functor $F$ with epimorphism $(-,
C)\longrightarrow F$, then for any subfunctor $F'$ of $F$, $F'=F$
provided that $F'(C)=F(C)$). Observe that a functor having finite
length is necessarily noetherian by an argument on its total length
(i.e., $l(F)=\sum_{C\in {\rm ind}(\mathcal{C})} l_R(F(C))$, where
$l_R$ denotes the length function on finitely generated
$R$-modules).

\vskip 10pt

The following result is essentially due to Auslander (compare
\cite{A}, Proposition 3.10).

\begin{lem}
Let $\mathcal{C}$ be a dualizing $R$-variety, $F\in
(\mathcal{C}^{\rm op}, R\mbox{\rm -mod})$. Then $F$ has finite
length if and only if $F$ is finitely presented and noetherian.
\end{lem}

\noindent{\bf Proof.}\quad Recall from \cite{AR2}, Corollary 3.3
that for a dualizing $R$-variety, functors having finite length are
finitely presented. So the ``only if'' follows.

 For the ``if'' part,
assume that $F$ is finitely presented and noetherian. Since $F$ is
finitely presented, by \cite{AR2}, p.324, we have the filtration of
subfunctors
\begin{align*}
0={\rm soc}_0(F) \subseteq {\rm soc}_1(F) \subseteq \cdots \subseteq
{\rm soc}_{i+1}(F) \subseteq \cdots
\end{align*}
where ${\rm soc}_1(F)$ is the socle of $F$, and in general ${\rm
soc}_{i+1}$ is the preimage of the socle of $F/{{\rm soc}_i(F)}$
under the canonical morphism $F \longrightarrow F/{{\rm soc}_i(F)}$.
Since $F$ is noetherian, we get ${\rm soc}_{i_0}{F}={\rm
soc}_{i_0+1}(F)$ for some $i_0$, and that is, the socle of $F/{{\rm
soc}_{i_0}(F)}$ is zero. However, by the dual of \cite{AR2},
Proposition 3.5, we know that for each nonzero finitely presented
functor $F$, the socle ${\rm soc}(F)$ is necessarily nonzero and
finitely generated semisimple. In particular, ${\rm soc}(F)$ has
finite length, and thus it is finitely presented. Note that ${\bf
fp}(C)\subseteq (\mathcal{C}^{\rm op}, R\mbox{-mod})$ is an abelian
subcategory, closed under extensions. Thus $F/{{\rm soc}_1(F)}$ is
finitely presented. Applying the above argument to $F/{{\rm
soc}_1(F)}$, we obtain that ${\rm soc}_2(F)$, as the extension
between the socles of two finitely presented functors, has finite
length. In general, one proves that $F/{\rm soc}_i(F)$ is finitely
presented and ${\rm soc}_{i+1}(F)$ has finite length for all $i$.
Hence ${\rm soc}(F/{\rm soc}_{i_0}(F))=0$ will imply  that $F/{\rm
soc}_{i_0}(F)=0$, i.e., $F={\rm soc}_{i_0}(F)$, which has finite
length. \hfill $\blacksquare$

\vskip 10pt

Let us consider the category $A\mbox{-GProj}$. Similar to Lemma
2.1(1),(2), we recall that $A\mbox{-GProj}\subseteq A\mbox{-Mod}$ is
closed under taking direct summands, kernels of epimorphisms and
extensions, and it is a Frobenius exact category with (relative)
projective-injective objects precisely contained in $A\mbox{-Proj}$.
Consider the stable category $A\underline{\mbox{-GProj}}$, which is
also triangulated and has arbitrary coproducts. Recall that in an
additive category $\mathcal{T}$ with arbitrary coproducts, an object
$T$ is said to be \emph{compact}, it the functor ${\rm
Hom}_\mathcal{T}(T, -)$ commutes with coproducts. Denote the full
subcategory of compact objects by $\mathcal{T}^c$. If we assume
further that $\mathcal{T}$ is triangulated, then $\mathcal{T}^c$ is
a thick triangulated subcategory. We say that $\mathcal{T}$ is a
\emph{compactly generated} \cite{Ne1,Ne}, if the subcategory
$\mathcal{T}^c$ is skeletally-small and for each object $X$,
$X\simeq 0$ provided that ${\rm Hom}_\mathcal{T}(T, X)=0$ for every
compact object $T$.

\vskip 5pt

 Note that in our situation, we always have an inclusion
$A\underline{\mbox{-Gproj}} \hookrightarrow
A\underline{\mbox{-GProj}}$, and in fact, we view it as
$A\underline{\mbox{-Gproj}}\subseteq
(A\underline{\mbox{-GProj}})^c$. Next lemma, probably known to
experts, states the converse in Gorenstein case. It is a special
case of \cite{C1}, Theorem 4.1 (compare \cite{Bel3}, Theorem 6.6).
One may note that in the artin case, the category
$A\underline{\mbox{-Gproj}}$  is idempotent-split. \vskip 15pt

\begin{lem}
Let $A$ be an Gorenstein artin algebra. Then the triangulated
category $A\underline{\mbox{\rm -GProj}}$ is compactly generated and
$A\underline{\mbox{\rm -Gproj}}\subseteq (A\underline{\mbox{\rm
-GProj}})^c$ is dense (i.e., surjective up to isomorphisms).
\end{lem}

\vskip 20pt

\subsection{Proof of Main Theorem:} Assume that $A$ is an artin
$R$-algebra. Set $\mathcal{C}=A\underline{\mbox{-Gproj}}$, by Lemma
2.2, $\mathcal{C}$ is a dualizing $R$-variety. For a finitely
generated Gorenstein-projective module $M$, we will denote by $(-,
M)$ the functor ${\rm Hom}_\mathcal{C}(-, M)$; for an arbitrary
module $X$, we denote by $(-, X)|_\mathcal{C}$ the restriction of
the functor $\underline{{\rm Hom}}_A(-,X)$ to $\mathcal{C}$.

 For the ``if" part,
we assume that each Gorenstein-projective module is a direct sum of
finitely generated ones. It suffices to show that the set ${\rm
ind}(\mathcal{C})$ is finite. For this end, assume that $M$ is a
finitely generated Gorenstein-projective module. We claim that the
functor $(-,M)$ is noetherian. In fact, given a subfunctor
$F\subseteq (-, M)$, first of all, we may find an epimorphism
\begin{align*}
\oplus_{i\in I} (-, M_i) \longrightarrow F,
\end{align*}
where each $M_i\in \mathcal{C}$ and $I$ is an index set. Compose
this epimorphism with the inclusion of $F$ into $(-, M)$, we get a
morphism from $\oplus_{i\in I}(-, M_i)$ to $(-,M)$. By the universal
property of coproducts and then by Yoneda's Lemma, we have, for each
$i$, a morphism $\theta_i: M_i \longrightarrow M$, such that $F$ is
the image of the morphism
\begin{align*}\sum_{i\in I}(-, \theta_i):
\oplus_{i\in I}(-, M_i)\longrightarrow (-,M).\end{align*} Note that
$\oplus_{i\in I}(-, M_i)\simeq (-, \oplus_{i\in
I}M_i)|_\mathcal{C}$, and the morphism above is also induced by the
morphism $\sum_{i\in I}\theta_i: \oplus_{i\in I}M_i \longrightarrow
M$. Form a triangle in $A\underline{\mbox{-GProj}}$
\begin{align*}
K[-1] \longrightarrow \oplus_{i\in I}M_i\stackrel{\sum_{i\in
I}\theta_i}\longrightarrow M \stackrel{\phi}\longrightarrow K.
\end{align*}
By assumption, we have a decomposition $K=\oplus_{j\in J}K_j$ where
each $K_j$ is finitely generated Gorenstein-projective. Since the
module $M$ is finitely generated, we infer that $\phi$ factors
through a finite sum $\oplus_{j\in J'}K_j$, where $J'\subseteq J$ is
a finite subset. In other words, $\phi$ is a direct sum of
\begin{align*}
M \stackrel{\phi'} \longrightarrow \oplus_{j\in J'}K_j \quad \mbox{
and } \quad 0 \longrightarrow \oplus_{j\in J\backslash J'}K_j.
\end{align*}
By the additivity of triangles, we deduce that there exists a
commutative diagram
\[\xymatrix@R=30pt{
\oplus_{i\in I} M_i \ar[r]^-{\sum_{i\in I} \theta_i} \ar[d] & M
\ar@{=>}[d]\\
M'\oplus (\oplus_{j\in J\backslash J'} K_j)[-1] \ar[r]^-{(\theta',
0)} & M }\]
 where the left side vertical map is an isomorphism, and
$M'$ and $\theta'$ are given by the triangle $(\oplus_{j\in
J'}K_j)[-1] \longrightarrow M' \stackrel{\theta'}\longrightarrow M
\stackrel{\phi'}\longrightarrow \oplus_{j\in J}K_j$. Note that
$M'\in \mathcal{C}$, and by the above diagram we infer that $F$ is
the image of the morphism $(-, \theta'): (-, M')\longrightarrow
(-,M)$, and thus $F$ is finitely-generated. This proves the claim.

By the claim, and by Lemma 2.3, we deduce that for each $M \in
\mathcal{C}$, the functor $(-, M)$ has finite length, in particular,
${\rm supp}((-, M))$ is finite. Assume that $\{S_i\}_{i=1}^n$ is a
complete list of pairwise nonisomorphic simple $A$-modules. Denote
by $f_i: X_i \longrightarrow S_i$ the right minimal $A\mbox{\rm
-Gproj}$-approximations. By Lemma 2.1(4), the module $M$ is a direct
summand of $M'$ and we have a finite chain of submodules of $M'$
with factors being among $X_i$'s. Then it is not hard to see that
${\rm supp}((-,M))\subseteq {\rm supp}((-,M'))\subseteq
\bigcup_{i=1}^n {\rm supp}((-, X_i))$ for every $M \in \mathcal{C}$.
Therefore we deduce that ${\rm ind}(\mathcal{C})= \bigcup_{i=1}^n
{\rm supp}((-, X_i))$, which is finite. \vskip 5pt

For the ``only if" part, assume that $A$ is a CM-finite Gorenstein
algebra. Then the set ${\rm ind}(\mathcal{C})$ is finite, say ${\rm
ind}(\mathcal{C})=\{G_1, G_2, \cdots, G_m\}$. Let $B={\rm
End}_\mathcal{C}(\oplus_{i=1}^m G_i)^{\rm op}$. Then $B$ is also an
artin $R$-algebra. Note that for each $C\in \mathcal{C}$, the
Hom-space ${\rm Hom}_\mathcal{C}(\oplus_{i=1}^mG_i, C)$ has a
natural left $B$-module structure, moreover, it is a finitely
generated projective $B$-module. In fact, we get an equivalence of
categories
\begin{align*}
\Phi={\rm Hom}_\mathcal{C}(\oplus_{i=1}^mG_i, -): \; \mathcal{C}
\longrightarrow B\mbox{-proj}.
\end{align*}
Then the equivalence above naturally induces the following
equivalences, still denoted by $\Phi$
\begin{align*}
\Phi:\; {\bf fp}(\mathcal{C}) \longrightarrow B\mbox{-mod},\quad
\Phi: \; (\mathcal{C}^{\rm op}, R\mbox{-Mod})\longrightarrow
B\mbox{-Mod}.
\end{align*}
In what follows, we will use these equivalences. By \cite{Ne}, p.169
(or \cite{C2}, Proposition 2.4), we know that the category ${\bf
fp}(\mathcal{C})$ is a Frobenius category. Therefore, via $\Phi$, we
get that $B$ is a self-injective algebra. Therefore by \cite{AF},
Theorem 31.9, we get that $B\mbox{-Mod}$ is also a Frobenius
category, and by \cite{AF}, p.319, every projective-injective
$B$-module is of form $\oplus_{i=1}^mQ_m^{(I_i)}$, where $\{Q_1,
Q_2, \cdots, Q_m\}$ is a complete set of indecomposable projective
$B$-modules such that $Q_i=\Phi(G_i)$, and each $I_i$ is some index
set, and $Q_i^{(I_i)}$ is the corresponding coproduct.

 Take $\{P_1, P_2, \cdots, P_n\}$ to be a complete set of
 pairwise nonisomorphic indecomposable projective $A$-modules. Let
 $G\in A\mbox{-GProj}$. We will show that $G$ is a direct sum of
 some copies of $G_i$'s and $P_j$'s. Then we are done. Consider the
 functor $(-, G)|_\mathcal{C}$, which is cohomological, and thus by
 \cite{C2}, Lemma 2.3 (or \cite{Ne}, p.258), we get
 ${\rm Ext}^1(F, (-, G)|_\mathcal{C})=0$ for each $F\in {\bf
 fp}(\mathcal{C})$, where the Ext group is taken in $(\mathcal{C}^{\rm op},
 R\mbox{-Mod})$. Via $\Phi$ and applying the Baer's criterion, we
 get that $(-, G)|_\mathcal{C}$ is an injective object, and thus by the
 above, we get an isomorphism of functors
 \begin{align*}
 \oplus_{i=1}^m (-, T_i)^{(I_i)} \longrightarrow (-,
 G)|_\mathcal{C},
 \end{align*}
where $I_i$ are some index sets. As in the first part of the proof,
we get a morphism $\theta:\oplus_{i=1}^mT_i^{(I_i)} \longrightarrow
T$ such that it induces the isomorphism above. Form the triangle in
$A\underline{\mbox{-GProj}}$
\begin{align*}
\oplus_{i=1}^m G_i^{(I_i)} \stackrel{\theta} \longrightarrow T
\longrightarrow X \longrightarrow (\oplus_{i=1}^m G_i^{(I_i)} )[1].
\end{align*}
For each $C\in \mathcal{C}$, applying the cohomological functor
${\rm Hom}_{A\underline{{\rm \mbox{-}GProj}}}(C,-)$ and by the
property of $\theta$, we obtain that
\begin{align*}{\rm
Hom}_{A\underline{{\rm\mbox{-}GProj}}}(C, X)=0,\quad \forall\; C\in
\mathcal{C}.
\end{align*}
By Lemma 2.4, the category $A\underline{\mbox{-GProj}}$ is generated
by $\mathcal{C}$, and thus $X\simeq 0$, and hence $\theta$ is an
isomorphism in the stable category $A\underline{\mbox{-GProj}}$.
Thus it is well-known (say, by \cite{CZ}, Lemma 1.1) that this will
force an isomorphism in the module category
\begin{align*}
\oplus_{i=1}^mG_i^{(I_i)}\oplus P\simeq G\oplus Q,
\end{align*} where $P$ and $Q$ are projective $A$-modules. Now by \cite{AF},
p.319, again, $P$ is a direct sum of copies of $P_j$'s. Hence the
combination of Azumaya's Theorem and Crawlay-J{\o}nsson-Warfield's
Theorem (\cite{AF}, Corollary 26.6) applies in our situation, and
thus we infer that $G$ is isomorphic to a direct sum of copies of
$G_i$'s and $P_j$'s. This completes the proof. \hfill $\blacksquare$

\vskip 10pt

\noindent {\bf Acknowledgement:}\; The author would like to thank
the referee very much for his/her helpful suggestions and comments.

\bibliography{}

\end{document}